\documentclass[runningheads]{llncs}

\usepackage{caption}
\usepackage{lmodern,textcomp}
\usepackage[T1]{fontenc}
\usepackage{graphicx}
\usepackage{float}
\usepackage[dvipsnames]{xcolor}
\usepackage{amsmath}
\usepackage{amssymb}
\usepackage[hidelinks]{hyperref}
\usepackage{stmaryrd}
\usepackage{booktabs}

\begin{document}
\title{Parareal with a physics-informed neural network as coarse propagator\thanks{This project has received funding from the European High-Performance Computing Joint Undertaking (JU) under grant agreement No 955701. The JU receives support from the European Union's Horizon 2020 research and innovation programme and Belgium, France, Germany, and Switzerland. This project also received funding from the German Federal Ministry of Education and Research (BMBF) grant 16HPC048.
This project has also received funding from the German Federal Ministry of Education and Research (BMBF) under grant 16ME0679K.}}
%
%
\author{Abdul Qadir Ibrahim\inst{1}\orcidID{0000-0003-3452-8500} \and
Sebastian G\"otschel\inst{1}\orcidID{0000-0003-0287-2120} \and
Daniel Ruprecht\inst{1}\orcidID{0000-0003-1904-2473}}
\authorrunning{A. Ibrahim et al.}
%
\institute{Chair Computational Mathematics, Institute of Mathematics, Hamburg University of Technology, Hamburg, Germany\\
\email{\{abdul.ibrahim,sebastian.goetschel,daniel.ruprecht\}@tuhh.de}}
\maketitle              
\begin{abstract} 
Parallel-in-time algorithms provide an additional layer of concurrency for the numerical integration of models based on time-depen\-dent differential equations. 
Methods like Parareal, which parallelize across multiple time steps, rely on a computationally cheap and coarse integrator to propagate information forward in time, while a parallelizable expensive fine propagator provides accuracy.
Typically, the coarse method is a numerical integrator using lower resolution, reduced order or a simplified model.
Our paper proposes to use a physics-informed neural network (PINN) instead.
We demonstrate for the Black-Scholes equation, a partial differential equation from computational finance, that Parareal with a PINN coarse propagator provides better speedup than a numerical coarse propagator.
Training and evaluating a neural network are both tasks whose computing patterns are well suited for GPUs. 
By contrast, mesh-based algorithms with their low computational intensity struggle to perform well.
We show that moving the coarse propagator PINN to a GPU while running the numerical fine propagator on the CPU further improves Parareal's single-node performance.
This suggests that integrating machine learning techniques into parallel-in-time integration methods and exploiting their differences in computing patterns might offer a way to better utilize heterogeneous architectures.

\end{abstract}

\section{Introduction}\label{sec:introduction}
Models based on differential equations are ubiquitous in science and engineering.
High-resolution requirements, often due to the multiscale nature of many problems, typically require that these models are run on high-performance computers to cope with memory demand and computational cost.
Spatial parallelization is already a widely used and effective approach to parallelize numerical algorithms for partial differential equations but, on its own, will not deliver enough concurrency for extreme-scale parallel architectures.
Parallel-in-time integration algorithms can help to increase the degree of parallelism in numerical models.
Combined space-time parallelization can improve speedup over spatial parallelization alone on hundreds of thousands of cores~\cite{SpeckEtAl2012}.

Parallel-in-time methods like Parareal~\cite{LionsEtAl2001}, PFASST~\cite{EmmettMinion2012} or MGRIT~\cite{FalgoutEtAl2014_MGRIT} rely on serial coarse level integrators to propagate information forward in time.
These coarse propagators constitute an unavoidable serial bottleneck which limits achievable speedup.
Therefore, the coarse-level integrators must be as fast as possible.
However, these methods are iterative and speedup will also decrease as the number of iterations goes up.
A coarse propagator that is too inaccurate, even when computationally cheap, will not provide good speedup because the number of required iterations will be too large.
Hence, a good coarse propagator needs to be at least somewhat accurate but also needs to run as fast as possible.
This trade-off suggests that using neural networks as coarse propagators could be promising: once trained, they are very fast to evaluate while still providing reasonable accuracy.
Furthermore, neural networks are well suited for running on GPUs whereas mesh-based discretizations are harder to run efficiently because of their lower computational intensity.
Therefore, algorithms featuring a combination of mesh-based components and neural network components would be well suited to run on heterogeneous systems combining CPUs and GPUs or other accelerators.

Our paper makes three novel contributions.
It (i) provides the first study of using a PINN as a coarse propagator in Parareal, (ii) shows that a PINN as a coarse propagator can accelerate Parareal convergence and improve speedup and (iii) illustrates that moving the PINN coarse propagator to a GPUs improves speedup further.
While we demonstrate our approach for the Black-Scholes equation, a model from computational finance, the idea is transferable to other types of partial differential equations where Parareal was shown to be effective.
We only investigate performance on a single node with one GPU.
Extending the approach to parallelize in time across multiple nodes and to work in combination with spatial parallelization left for future work.
The code used to generate the results shown in this paper is freely available~\cite{pinn_blackscholes}.

\section{Related Work}\label{sec:related-work}
Using machine learning (ML) to solve differential equations has become an active field of research. 
Some papers aim to entirely replace the numerical solver by neural networks~\cite{RanadeEtAl2021,StenderEtAl2022}.
Physics-informed neural networks (PINNs)~\cite{raissi2017physics}, which use the residual of a partial differential equation (PDE) as well as boundary- and initial conditions in the loss function, are used in many applications.
This includes a demonstration for the Black Scholes equation~\eqref{eq:black_scholes_equations}, showing that a PINN is capable of accurately pricing a range of options with complex payoffs, and is significantly faster than traditional numerical methods~\cite{sirignano2018dgm}.
However, solving differential equations with ML alone generally does not provide the high accuracy that can be achieved by numerical solvers.
This has led to a range of ideas where ML is used as an ingredient of classical numerical methods instead and not as a replacement~\cite{HuangEtAl2022}.

Specific to parallel-in-time integration methods, there are two research directions aiming to connect them with machine learning.
On the one hand, there are attempts to use ML techniques to improve parallel-in-time algorithms.
Our paper falls into this category.
Using a neural network as coarse propagator for Parareal has been studied in two previous papers.
Yalla and Enquist~\cite{YallaEnquist2018} were the first to explore this approach.
They use a neural network with one hidden layer of size 1000 and demonstrate for a high dimensional oscillator that it helps Parareal converge faster compared to a numerical coarse propagator.
However, no runtimes or speedups are reported.
Agboh et al.~\cite{AgbohEtAl2020} use a feed-forward deep neural network as a coarse propagator to integrate an ordinary differential equation modeling responses to a robot arm pushing multiple objects.
They also observe that the trained coarse propagator improves Parareal convergence compared to a simplified analytical coarse model.
Nguyen and Tsai~\cite{NguyenTsai2023} do not fully replace the numerical coarse propagator but use supervised learning to enhance its accuracy for wave propagation modeling. 
They observe that this enhances stability and accuracy of Parareal, provided the training data contains sufficiently representative examples. 
Gorynina et al.~\cite{GoryninaEtAl2022} study the use of a machine-learned spectral neighbor analysis potential in molecular dynamics simulations with Parareal.

A few papers go the opposite way and adopt ideas from parallel-in-time integration methods to parallelize and accelerate the process of training deep neural networks.
G\"unther et al.~\cite{GuentherEtAl2020} use a nonlinear multi-grid method to improve the training process of a deep residual network.
They use MGRIT, a multi-level generalization of Parareal, to obtain layer-parallel training on CPUs, reporting a speedup of up to 8.5 on 128 cores. 
Kirby et al.~\cite{KirbyEtAl2020} extend their approach to multiple GPUs, obtaining further performance gains.
In a similar way, Meng et al.~\cite{MengEtAl2020} use Parareal to generate starting values for a series of PINNs to help with the training process. 
Motivated by the observation that it becomes expensive to train PINNs that integrate over long time intervals, they concatenate multiple short-time PINNs instead.
They use a cheap numerical coarse propagator and a Parareal iteration to connect these PINNs with each PINN inheriting the parameters from its predecessor.
While they mention the possibility of using a PINN as coarse propagator, they do not pursue this idea further in their paper.
Lorin~\cite{Lorin2020} derives a parallel-in-time variant of neural ODEs to improve training of deep Residual Neural Networks. 
Finally, Lee et al.~\cite{LeeEtAl2022} use a Parareal-like procedure to train deep neural networks across multiple GPUs.

\section{Algorithms and Benchmark Problem}\label{sec:algorithms-and-benchmark-problem}
The Black-Scholes equation is a widely used model to price options in financial markets~\cite{black1973pricing}.
It is based on the assumption that the price of an asset follows a geometric Brownian motion, so that the log-returns of the asset are normally distributed.
Closed form solutions exist for the price of a European call or put option~\cite{kumar2012analytical}, but not for more complex options such as American options or options with multiple underlying assets.
To be able to compute numerical errors, we thus focus on the European call option, a financial derivative that gives the buyer the right, but not the obligation, to buy an underlying asset at a predetermined price (the strike price) on or before the expiration date.
The price $V$ of the option can be modeled by
\begin{equation}
   f(V) = \frac{\partial V}{\partial t}(S, t) + \frac{1}{2} \sigma^2 S^2 \frac{\partial^2 V}{\partial S^2}(S, t) + rS \frac{\partial V}{\partial S}(S, t) - rV(S, t) = 0,
    \label{eq:black_scholes_equations}
\end{equation}
where $S$ denotes the current value of the underlying asset, $t$ is time, $r$ denotes the no-risk interest rate (for example saving rates in a bank) and $\sigma$ denotes the volatility of the underlying asset.
To fully determine the solution to~\eqref{eq:black_scholes_equations}, we impose a final state at expiry time $t=T$ and two boundary conditions with respect to $S$, motivated by the behaviour of the option at $S = 0$ and  as $S \rightarrow \infty$.
For the call option, the expiry time condition is
\begin{equation}
    \label{eq:boundary_conditions_a}
    V(T, S) = \max(S - K, 0) \ \text{for all} \ S.
\end{equation}
If the underlying asset becomes worthless, then it will remain worthless, so the option will also be worthless.
Thus, 
\begin{equation}
    \label{eq:boundary_conditions_b}
    V(t, 0) = 0 \ \text{for all} \ t.
\end{equation}
On the other hand, if $S$ becomes very large, then the option will almost certainly
be exercised, and the exercise price is negligible compared to $S$.
Thus, the option will have essentially the same value as the underlying asset itself and
\begin{equation}
    \label{eq:boundary_conditions_c}
    V(t, S) \sim 0 \ \text{as} \ S \rightarrow \infty, \ \text{for fixed} \ t.
\end{equation}
For the European call option, we select an interval of $t =0$ and $T=1$ and an artificial bound for the asset of $S = 5000$\texteuro.

\subsection{Parareal}\label{subsec:parareal}
Parareal is an iterative algorithm to solve an initial value problem of the form
\begin{equation}
	 V'(t) = \phi(V(t)), \ t \in [0,T], \ V(0) = V_0,
    \label{eq:ivp}
\end{equation}
where in our case the right hand side function $\phi$ stems from the discretization of the spatial derivatives in~\eqref{eq:black_scholes_equations}.
Note that the coefficients in~\eqref{eq:black_scholes_equations} do not depend on time, so we can restrict our exposition to the autonomous case.
Decompose the time domain $[0,T]$ into $N$ time-slices $[T^n, T^{n+1}]$, $n=0, \ldots, N-1$.
Denote as $\mathcal{F}$ a numerical time stepping algorithm with constant step size $\delta t$ and high accuracy and as
\begin{equation}
	\label{eq:serial-fine}
	V_{n+1} = \mathcal{F}(V_n)
\end{equation}
the result of integrating from some initial value $V_n$ at the start time $T^n$ of a time slice until the end time $T^{n+1}$.
Classical time stepping corresponds to evaluating~\eqref{eq:serial-fine} for $n=0, \ldots, N-1$ in serial.
Parareal replaces this serial procedure with the iteration
\begin{equation}
	\label{eq:parareal}
	V^{k+1}_{n+1} = \mathcal{G}(V^{k+1}_n) + \mathcal{F}(V^k_n) - \mathcal{G}(V^k_n)
\end{equation}
where $k=1, \ldots, K$ counts the iterations.
The key in~\eqref{eq:parareal} is that the computationally expensive evaluation of $\mathcal{F}$ can be parallelized across all $N$ time slices.
Here, we always assume that $P=N$ many processes are used and each process holds a single time slice.
A visualization of the Parareal workflow as well as pseudocode can be found in the literature~\cite{Ruprecht2017_lncs}.
As $k \to N$, $V^{k}_n$ converges to the same solution generated by serial evaluation of~\eqref{eq:serial-fine}.
However, to achieve speedup, we require convergence in $K \ll N$ iterations.
An upper bound for speedup achievable with Parareal using $P$ processors to integrate over $N=P$ time slices is given by
\begin{equation}
	\label{eq:speedup}
	s_{\text{bound}}(P) = \frac{1}{ \left(1 + \frac{K}{P} \right) \frac{c_{\text{c}}}{c_{\text{f}}} + \frac{K}{P} }
\end{equation}
where $K$ is the number of iterations, $c_{\text{c}}$ the runtime of $\mathcal{G}$ and $c_{\text{f}}$ the runtime of $\mathcal{F}$~\cite{Ruprecht2017_lncs}.
Since~\eqref{eq:speedup} neglects overhead and communication, it is an upper bound on achievable speedups and measured speedups will be lower.

\subsection{Numerical solution of the Black-Scholes equation}\label{subsec:mesh-based-discretization}
We approximate the spatial derivatives in~\eqref{eq:black_scholes_equations} by second order centered finite differences on an equidistant mesh
\begin{equation}
	0 = S_0 < S_1 < \ldots < S_N = L
\end{equation}
with $S_{i+1} - S_i = \Delta S$ for $i=0, \ldots, N-1$.
For the inner nodes, we obtain the semi-discrete initial value problem
\begin{equation}
	V^{'}_j(t) = -\frac{1}{2} \sigma^2 S_j^2 \frac{V_{j+1} - 2 V_j + V_{j-1}}{\Delta S^2} - r S_j \frac{V_{j+1} - V_{j-1}}{2 \Delta S} + r V_j
	\label{eq:semi-discrete-equation}
\end{equation}
with $j=1, \ldots, $.
This is complemented by the boundary condition $V_0 = 0$ for a zero asset value.
We also impose the asymptotic boundary condition~\eqref{eq:boundary_conditions_c} at finite distance $L$ so that $V_N = 0$.
In time, we use a second order Crank-Nicolson method for $\mathcal{F}$ and a first order implicit Euler method as numerical $\mathcal{G}$.
Since we have a final condition instead of an initial condition, we start at time $T = 1$ and solve the problem backwards.
We use $200$ steps for the fine method and $100$ steps for the coarse.

\subsection{Physics Informed Neural Network (PINN)}\label{subsec:physics-informed-neural-network-(pinn)}
\begin{figure}[ht]
    \centering
    \centering
    \includegraphics[scale=1]{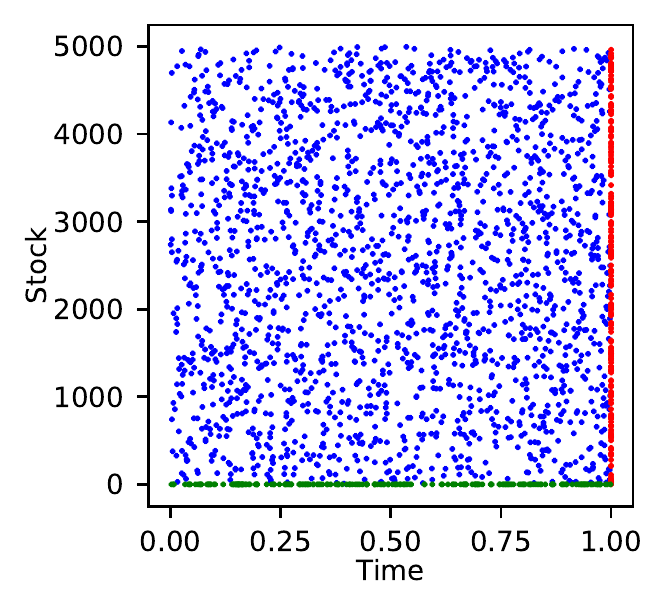}
    \caption{Subset of the randomly generated collocation nodes. The solution is forced to satisfy the PDE at the inner nodes by minimizing the PDE residual, to satisfy the boundary condition at the green nodes via the boundary loss and the expiration condition at the red nodes via the expiration loss.}
    \label{fig:collocation}
\end{figure}

The PINN we use as coarse propagator gets a time slice $[t_\text{start}, t_\text{end}] \subset  [0, T]$, the asset price $V$ at $t_\text{start}$ and stock values $S$, and  outputs the predicted state of the asset price $\tilde{V}$ at $t_\text{end}$.
To train it, we define three sets of collocation points in time and stock price: $(S_i, t_i), i=1,\dots N_f$ in the interior of the space-time domain for evaluating the residual $f(V)$ of the Black-Scholes eqation \eqref{eq:black_scholes_equations}, 
 $(S_i, t_i), i=1,\dots N_b$ collocation points on the boundary to evaluate \eqref{eq:boundary_conditions_a}, and $S_i, i=1,\dots N_\text{exp}$ for the final state conditions \eqref{eq:boundary_conditions_b}, \eqref{eq:boundary_conditions_c}. 
The loss function to be minimized is given by
\begin{equation}
\text{MSE}_{\text{total}} = \text{MSE}_f + \text{MSE}_{\exp} + \text{MSE}_b,
\label{eq:bs_totalloss}
\end{equation}
consisting of a term to minimize the PDE residual $f(V)$
\begin{align}
\text{MSE}_{\text{f}} &=  \frac{1}{N_f}\sum_{i=1}^{N_f} |f(\tilde{V}(t_i, S_i))|^2,
\label{eq:loss_pde_call}
\end{align}
the boundary loss term
\begin{equation}
\text{MSE}_{\text{b}} = \frac{1}{N_b}\sum_{i=1}^{N_b}\left|\tilde{V}(t_i, S_i) - V(t_i, S_i)\right|^2,
\label{eq:loss_bc_call}
\end{equation}
and the loss at expiration
\begin{equation}
\text{MSE}_{\exp} = \frac{1}{N_{\exp}}\sum_{i=1}^{N_{\exp}}\left|\tilde{V}(T, S_i) - \max(S_i - K, 0)\right|^2,
\label{eq:loss_fs_call}
\end{equation}
For our setup, we randomly generate $N_f=100,000$ collocation points within the domain $[0,5000] \times [0,1]$, $N_b= 10,000$ collocation points at the boundary $[0,1]$ and $N_\text{exp}=10,000$ collocation points to sample the expiration condition over $[0,5000]$.
The derivatives that are required to compute the PDE loss are calculated by automatic differentiation~\cite{baydin2018automatic}.
We compute the PDE residual~\eqref{eq:loss_pde_call} over the points inside the domain, the boundary condition loss~\eqref{eq:loss_bc_call} over the spatial boundary and the expiration loss~\eqref{eq:loss_fs_call} over the end points.
The sum of the three forms the total loss function~\eqref{eq:bs_totalloss}.
Figure~\ref{fig:collocation} shows a subset of the generated collocation points to illustrate the approach.

\begin{figure}[th]
    \centering
    \includegraphics[width=\textwidth]{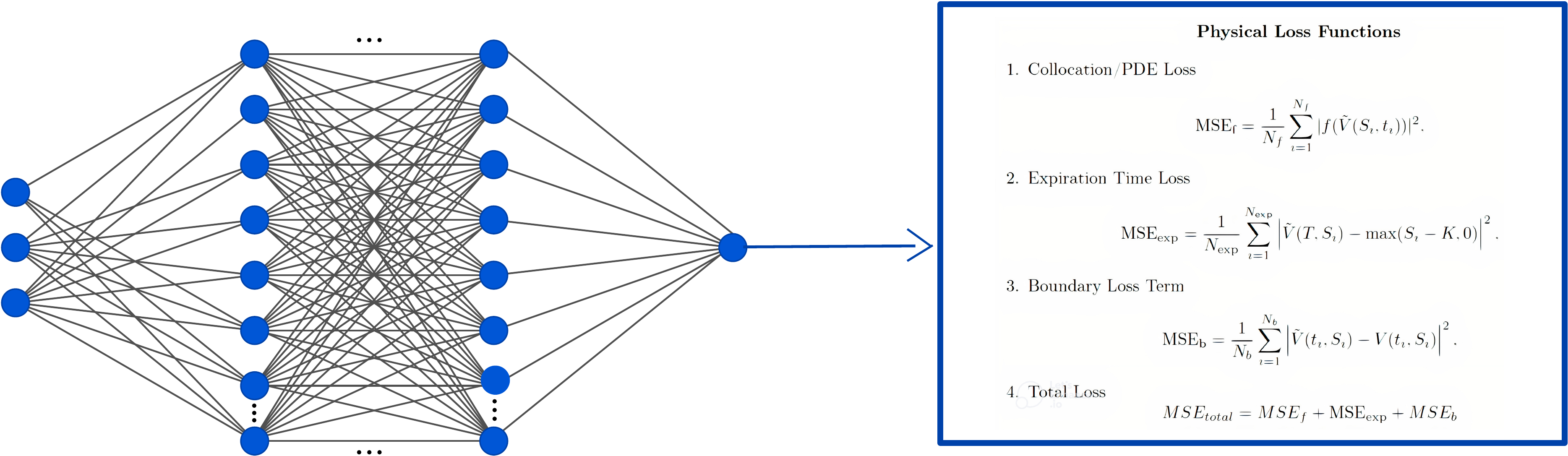}
    \caption{Structure of the PINN. The network takes the time $t_{\text{start}}, t_{\text{end}}$, asset values $V$ and stock values $S$ as input and returns the predicted asset
    values $\tilde{V}$ at $t_\text{end}$. The loss function encodes the PDE, the expiration condition and the boundary conditions. Figure produced using \url{https://alexlenail.me/NN-SVG/index.html}.}
    \label{fig:pinn_architecture}
\end{figure}

The neural network consists of 10 fully connected layers with 50 neurons in each and was implemented using Pytorch~\cite{paszke2019pytorch}. 
Figure~\ref{fig:pinn_architecture} shows the principle of a PINN but for a smaller network for the sake of readability.
Every linear layer, excluding the output layer, is followed by the ReLU activation function.
The weights for the neural network are initialized using Kaiming~\cite{he2015delving}.
We focus here on a proof-of-concept and have not undertaken a systematic effort to optimize the network architecture
but this would be an interesting avenue for future research.

We used the Adam optimizer~\cite{KingmaBa2014} with a learning rate of $10^{-2}$ for the initial round of training for 5000 epochs, followed by a second round of training with a learning rate of $10^{-3}$ for 800 epochs.
The training data (collocation points) was shuffled during every epoch to prevent the model from improving predictions based on data order rather than the underlying patterns in the data.
Table~\ref{tab:model_convergence_first} shows the behavior of the three loss function terms.
The total training time for this model was around 30 minutes.
\begin{table}[th]
    \renewcommand{\arraystretch}{1.5}
      \setlength\tabcolsep{4mm}
    \centering
    \begin{tabular}{cccc}
    	\toprule
         \textbf{ Epoch} & \textbf{Expiration}            & \textbf{Boundary}        & \textbf{Residual}                  \\ \midrule
         $0$  & $9.21 \times 10^{2}$  & $9.21 \times 10^{2}$        & $7.33 \times 10^{3}$             \\ 
         $2000$ & $5.58 \times 10^{-1}$ & $3. 45 \times 10^{-2}$        & $2.50 \times 10^{-2}$ \\ 
         $4000$ & $4.11 \times 10^{-2}$  & $2.34 \times 10^{-2}$        & $5.00 \times 10^{-3}$               \\ 
         $5000$ & $5.92 \times 10^{-1}$   & $1.34 \times 10^{-2}$ & $4.22 \times 10^{-3}$                  \\  \midrule
         $5300$ & $4.19 \times 10^{-2}$ & $3.22 \times 10^{-3}$          & $1.94 \times 10^{-4}$  \\ 
         $5500$ & $6.46 \times 10^{-4}$  & $1.96 \times 10^{-4}$         & $5.73 \times 10^{-5}$  \\ 
         $5800$ & $2.92 \times 10^{-5}$      & $1.14 \times 10^{-5} $  & $3.19 \times 10^{-4}$      \\ \bottomrule
    \end{tabular}
    \caption{Evolution of the loss function during network training.
    The three columns show the MSE for the three terms of the loss function related to the end condition~\eqref{eq:boundary_conditions_a}, boundary conditions~\eqref{eq:boundary_conditions_b} and~\eqref{eq:boundary_conditions_c} and residual~\eqref{eq:black_scholes_equations}.
    After 5000 epochs with training rate $10^{-2}$, another 800 epochs of training with a reduced training rate of $10^{-3}$ were performed.
}
    \label{tab:model_convergence_first}
\end{table}

\section{Results}\label{sec:results}
The numerical experiments were conducted on OpenSUSE Leap 15.4 running an Intel Core 24 x 12th Gen Intel i9-12900K with
a base clock speed of 3.2 GHz and a maximum turbo frequency of 5.2 GHz, with 62.6 GiB of RAM and an NVIDIA
GeForce RTX 3060/PCIe/SSE2 GPU. Implementations were done using Python 3.10,
pytorch1.13.1+cu117, mpi4py3.1.4, as well as numba0.55.1 for the GPU runs.

\begin{figure}[ht]
    \centering
        \centering
        \includegraphics[scale=.95]{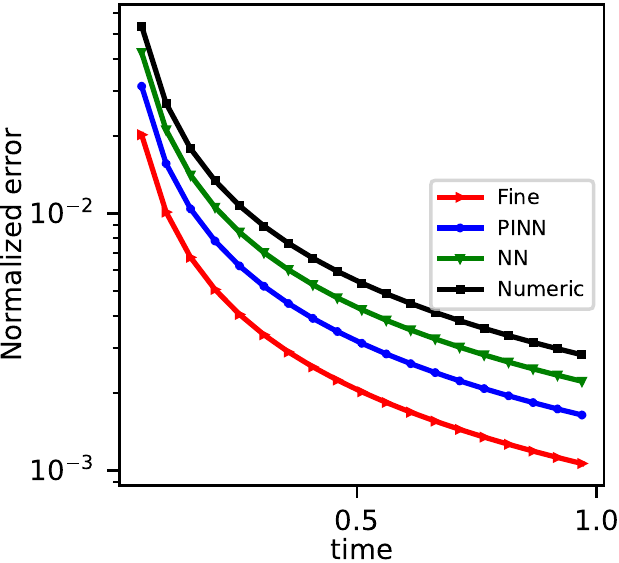}\hfill
        \includegraphics[scale=.95]{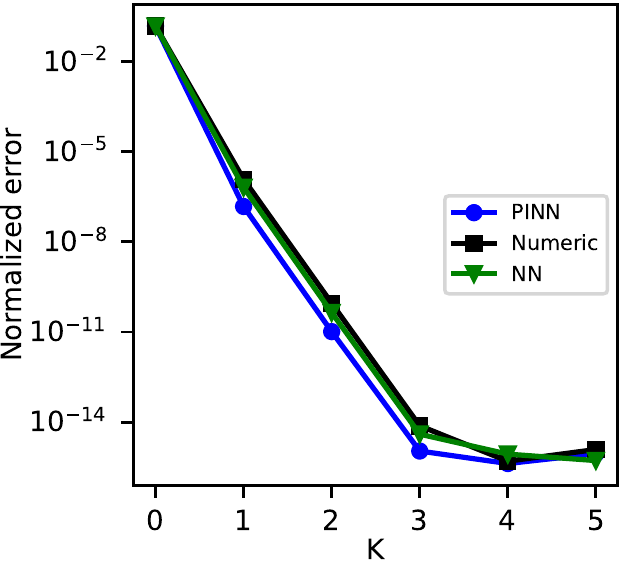}
        \caption{Normalized $\ell_2$-error over time of coarse and fine propagator against the analytical solution (left). Normalized $\ell_2$-error against the serial fine solution versus number of iterations for three different variants of Parareal (right). The black line (squares) is Parareal with a numerical coarse propagator, the green line (diamonds) is Parareal with a neural network as coarse propagator that is trained only on data while the blue line (circles) is Parareal with a PINN as coarse propagator that also uses the terms of the differential equation in the loss function. Parareal uses $P=16$ time slices in all cases.}
        \label{fig:parareal_pinn_nn}
\end{figure}

\paragraph{Parareal convergence.}
Figure~\ref{fig:parareal_pinn_nn} shows the normalized $\ell_2$ error for the serial fine, numerical coarse and PINN-coarse propagator over time (left).
As expected, the fine propagator is the most accurate with an $\ell_2$ error of around $10^{-3}$ at the end of the simulation.
The numerical coarse propagator is noticeably less accurate.
The PINN coarse propagator is more accurate than the numerical coarse propagator but also does not reach the accuracy of the fine.
To illustrate the importance of encoding the differential equation in the loss function, we also show a neural network (NN) trained only on data produced with the fine propagator but without the terms encoding the differential equation.
The neural network without PDE residual is somewhat more accurate than the numerical coarse method but not as good as the PINN. Note that the PINN used here does not need numerically generated trajectories as training data, as the loss function~\eqref{eq:bs_totalloss} only consists of PDE residual, boundary and expiration conditions and does not include a data mismatch term.

Figure~\ref{fig:parareal_pinn_nn} (right) shows the normalized $\ell_2$ error of Parareal against the number of iterations.
For all three coarse propagators, numerical, NN and PINN, Parareal converges very quickly.
Although PINN and NN are slightly more accurate than the numerical coarse propagator, the impact on convergence is small.
After one iteration, the iteration error of Parareal is smaller than the discretization error of the fine method.
After $K=3$ iterations, Parareal has reproduced the fine solution up to round-off error.
Below, we report runtimes and speedup for $K=3$.
With only a single iteration, the $K / P$ term in~\eqref{eq:speedup} is less important and reducing the runtime of the coarse propagator increases overall speedup even more.
Therefore, the case with $K=3$  is the case where switching to the coarse propagator will yield less improvement.

\begin{figure}[ht]
    \centering
    \centering
    \includegraphics[scale=0.95]{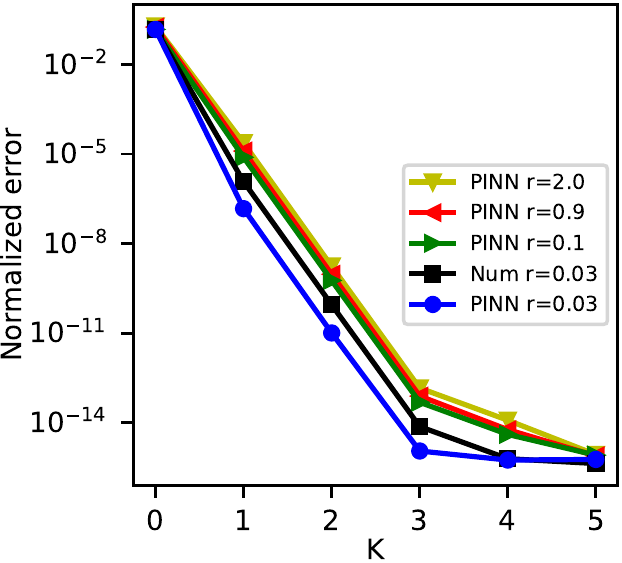}\hfill
    \includegraphics[scale=0.95]{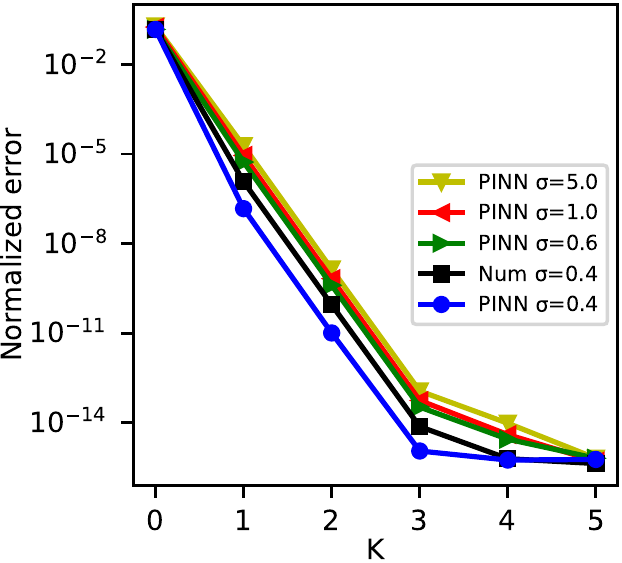}
    \caption{Convergence of Parareal for different interest rates $r$ (left) and volatilities $\sigma$ (right). In all cases, the coarse propagator is the PINN trained for values of $r=0.03$ and $\sigma=0.4$. Even for parameter values more than ten times larger than the ones for which the PINN was trained, Parareal requires only one additional iteration to converge to within machine precision of the fine integrator.}
    \label{fig:pinn-generalization}
\end{figure}

\paragraph{Generalization.}
Figure~\ref{fig:pinn-generalization} shows how Parareal with a PINN coarse propagator converges if applied to~\eqref{eq:black_scholes_equations} with parameters different from those for which the PINN was trained.
As parameters become increasingly different from the training values, the coarse propagator will become less accurate.
However, if Parareal converges, it will produce the correct solution since the numerical fine propagator always uses the correct parameters. 
The combination of Parareal + PINN generalizes fairly well.
Even for parameters more than ten times larger than the training values it only requires one additional iteration to converge.
While the additional iteration will somewhat reduce achievable speedup as given by~\eqref{eq:speedup}, the performance results presented below should not be overly sensitive to changes in the model parameters.

\paragraph{Parareal runtimes and speedup.}
Reported runtimes are measured using the \texttt{time} command in Linux and include the time required for setup,
computation and data movement.
Table~\ref{tab:numeric_vs_pinn} shows the runtime in milliseconds of Parareal using $P=16$ cores for four different coarse propagator configurations.
Shown are averages over five runs as well as the standard deviation.
Replacing the numerical coarse propagator with a PINN on a CPU reduces Parareal execution time by a factor of $2.4$, increasing to $2.9$ if the PINN is run on a GPU.
For the numerical coarse propagator, using the GPU offers no performance gain because the resolution and thus computational intensity is not high enough.
The much faster coarse propagator provided by the PINN significantly reduces the serial bottleneck in Parareal and will, as demonstrated below, yield a marked improvement in speedup.
\begin{table}[th]
    \renewcommand{\arraystretch}{1.5}
    \setlength\tabcolsep{2mm}
    \centering
    \begin{tabular}{lccc}
        \toprule
        \textbf{}                 & \textbf{Numerical}        & \textbf{PINN}              & Speedup over CPU-Numerical \\ \midrule
        CPU & $3.48  \pm 0.056 $  & $1.47  \pm 0.073$ & $2.4$\\
        GPU & $3.99  \pm 0.651 $  &$1.21  \pm 0.041$  & $2.9$ \\ 
        Speedup &   --                     &  $1.21$                     & \\ \bottomrule
    \end{tabular}
    \caption{Runtime $c_{\text{c}}$ in milliseconds of the coarse propagator $\mathcal{C}$ averaged over five runs plus/minus standard deviation.}
    \label{tab:numeric_vs_pinn}
\end{table}

Table~\ref{tab:pinngpu_vs_pinncpu} shows runtimes for the full Parareal iteration averaged over five runs.
The fastest configuration is the one that runs the numerical fine propagator on the CPU and the PINN coarse propagator on the GPU.
Executing both fine and coarse propagator on the CPU takes about a factor of three longer.
Importantly, moving both to the GPU, while somewhat faster than running all on the CPU, is slower than the mixed version by a factor of about two.
The full GPU variant will eventually be faster if the resolution of the fine and coarse propagator are both extremely high.
However, the current resolution already produces an error of around $10^{-3}$ which will be sufficient in most situations.
This illustrates how a combination of numerical method and PINN within Parareal can not only improve performance due to the lower cost of the PINN but also help to better utilize a node that features both CPUs and GPUs or even neural network accelerators.
Thus, the different computing patters in finite difference numerical methods and neural networks can be turned into an advantage.
\begin{table}[th]
    \renewcommand{\arraystretch}{1.5}
     \setlength\tabcolsep{2mm}
    \centering
    \begin{tabular}{lcc}
        \toprule
        \textbf{}                 & \textbf{CPU-Coarse}        & \textbf{GPU-Coarse}              \\ \midrule
        \textbf{CPU-Fine} & $128.48  \pm 0.715$  & $ 41.241970 \pm 0.334$\\ 
        \textbf{GPU-Fine} & $83.2545  \pm 0.356$ & $87.45234 \pm 0.253$       \\ \bottomrule
    \end{tabular}
    \caption{Runtimes in milliseconds for Parareal averaged over five runs plus/minus standard deviation.}
    \label{tab:pinngpu_vs_pinncpu}
\end{table}

\begin{figure}[th]
    \centering
    \centering
    \includegraphics[scale=0.95]{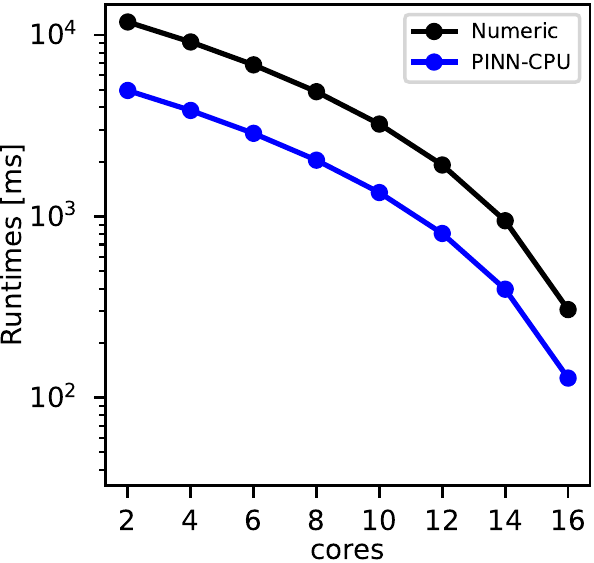}\hfill
    \includegraphics[scale=0.95]{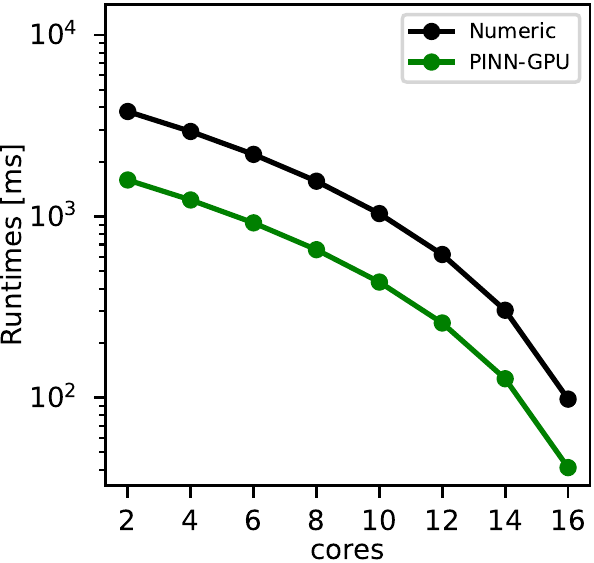}
    \caption{Runtimes in milliseconds for Parareal (dots) and the serial numerical fine propagator (horizontal lines)
        on a CPU (left) and GPU (right)}
    \label{fig:runtime-pinn}
\end{figure}

Figure~\ref{fig:runtime-pinn} shows runtimes for Parareal with both a PINN and numerical coarse propagator on a CPU (left) and GPU (right) against the number of cores/time slices $P$.
The numerical fine propagator is always run on the CPU.
In both cases, runtimes decrease at a similar rate as the number of time slices/cores $P$ increases.
The numerical coarse propagator is consistently slower than the PINN and the gap is similar on the CPU and GPU.
\begin{figure}[th]
    \centering
    \centering
    \includegraphics[scale=.95]{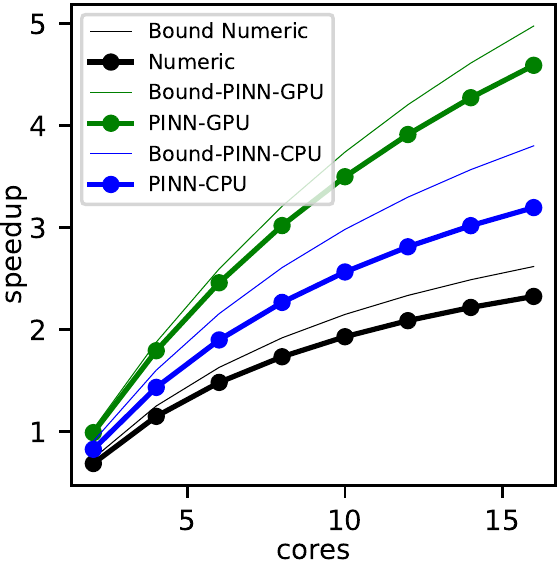}\hfill
    \includegraphics[scale=.95]{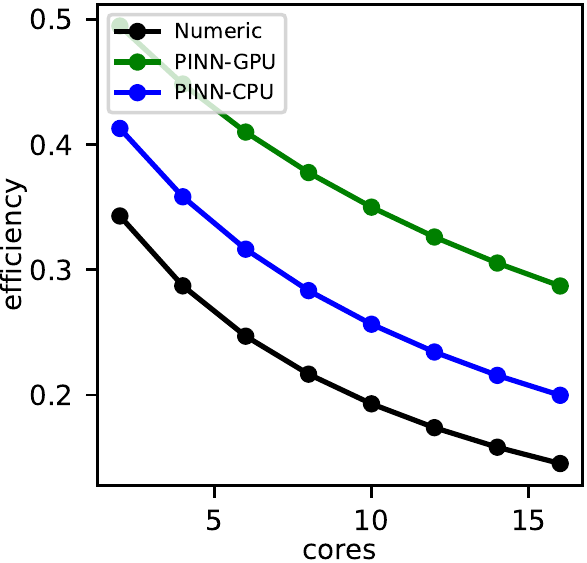}
    \caption{Speedup (left) and parallel efficiency (right) of Parareal over the serial numerical fine propagator on a CPU.
    Because the PINN-GPU coarse propagator is faster, it reduces the serial bottleneck of Parareal and allows for better speedup and parallel efficiency.}
    \label{fig:speedup_efficiency}
\end{figure}
Finally, Figure~\ref{fig:speedup_efficiency} shows the speedup (left) and parallel efficiency (right) for Parareal with a numerical, PINN-CPU and PINN-GPU coarse propagator.
The speedup bounds~\eqref{eq:speedup} are shown as lines.
Moving from a numerical coarse propagator to a PINN and moving the PINN from the CPU to a GPU each improves speedup significantly.
For the numerical coarse propagator, Parareal achieves a speedup of around $S(16) \approx 2$.
Replacing the numerical integrator with a PINN  improves speedup to $S(16) \approx 3$.
Running this PINN on a GPU again improves speedup to $S(16) \approx 4.5$, more than double what we achieved with the numerical coarse propagator on a CPU.
The improvements in speedup translate into increased parallel efficiency, which improves from around $30 \%$ for the numerical coarse propagator to around $60 \%$ for the PINN-GPU coarse method.
For smaller numbers of processors, the gains in speedup are less pronounced, because the $K/P$ term in~\eqref{eq:speedup} is more dominant.
But gains in parallel efficiency are fairly consistent from $P=2$ cores to $P=16$ cores.
In summary, this demonstrates that replacing a CPU-run numerical coarse propagator with a GPU-run PINN can greatly improve the performance of Parareal by minimizing the serial bottleneck from the coarse propagator.

\section{Discussion}\label{sec:discussion}
Parareal is a parallel-in-time method that iterates between a cheap coarse and a parallel expensive fine integrator.
To maintain causality, the coarse propagator needs to run in serial and therefore reflects a bottleneck that limits achievable speedup.
Mostly, coarse propagators are similar to fine propagators and build using numerical methods but with lower order, lower resolution or, in some cases, models of reduced complexity.
We investigate the use of a physics-informed neural network (PINN) instead.
The PINN is shown to be slightly more accurate than a numerical coarse propagator but a factor of three faster.
Using it does not affect convergence speed of Parareal but greatly reduces the serial bottleneck from the coarse propagator.

We show that, on a single node with one GPU, a combination of a numerical fine propagator run on a CPU with a PINN coarse propagator run on a GPU provides more than twice the speedup than vanilla Parareal using a numerical coarse propagator run on the CPU.
Also, we demonstrate that moving both fine and coarse propagator to the GPU is slower than moving just the PINN coarse method to the GPU and keeping the numerical fine method on the CPU.
The reason is that unless the resolution of the fine propagator is extremely high, its low computational intensity means there is little gain from computing on a GPU and so overheads from data movement are dominant.
By contrast, evaluating PINNs is well suited for GPU computation.
Our results demonstrate that using PINNs to build coarse level models for parallel-in-time methods is a promising approach to reduce the serial bottleneck imposed by causality. 
They also suggest that parallel-in-time methods featuring a combination of numerical algorithms and neural networks might be useful to better utilize heterogeneous systems.

\bibliographystyle{splncs04}
\bibliography{pint_mod}
\end{document}